\documentclass[12pt,preprint,nofootinbib,nobibnotes,nonumbib]{article}
\usepackage{amssymb}
\usepackage{amsmath}
\usepackage{hyperref}

\begin{document}

\title{Extending Representations of Dense Subalgebras of 
$C^\star$-Algebras, and Spectral Invariance}

\author{Larry B. Schweitzer}

\date{April 6, 2018}

\maketitle

\begin{abstract}
Let $A$ be a dense Fr\'echet subalgebra
of the $C^\star$-algebra of compact operators $\cal K$
on a seprable Hilbert space.
Assume that $A$ is spectral invariant in $\cal K$.  
We show that every algebraically
cyclic subrepresentation of a topologically irreducible
representation 
of $A$ is  contained in a $\cal K$-module.


\end{abstract}

\noindent
See [Sch, 2015] for definitions and introduction.

\noindent {\bf Theorem 1.1.}\ 
{\it Let $A$ be a dense Fr\'echet subalgebra of the $C^\star$-algebra of
compact operators $\cal K$ on an infinite dimensional Hilbert space.
Assume that $A$ is spectral invariant in $\cal K$.
If $E$ is any topologically irreducible Banach $A$-module,
then any algebraically cyclic $A$-submodule is contained in a
$\cal K$-module.}

\noindent {\bf Proof:}
Let $E$ be a topologically irreducible Banach $A$-module.
Let $F$ be an algebraically cyclic $A$-submodule of $E$ with
nonzero cyclic vector $\xi \in F$ and $b \in A$ satisfying $b\xi=\xi$.
Since $b\xi=\xi$, $1 \in Sp_A(b)$,
and since $A$ is spectral invariant in $\cal K$,
$1$ is an isolated point of $Sp_A(b)$.
By a standard holomorphic functional
calculus argument,
we can find an idempotent $e$ in the closed
subalgebra of $A$ generated by $b$, which satisfies $e\xi= \xi$.
(See for example [Palm, 1994], Proposition 3.4.1,
 [Sch, 1992], Lemma 1.2.)
Since $e$  is an idempotent in $\cal K$, $e$ has finite rank.

Since the sublagebra $e {\cal K} e$ of $\cal K$ is a full
finite dimensional matrix algebra,
and $e Ae$ is dense in it, the latter must
also be a full matrix algebra.  Hence we can find
rank 1 (minimal) idempotents $e_1, \dots e_n \in A$
which satisfy $e_i e_j = 0$ if
$i \ne j$, and $e_1 + e_2 + \cdots e_n = e$.
The left ideal $Ae$ of $A$ is the direct sum of
minimal left ideals $Ae_1 + Ae_2 + \cdots + Ae_n$.

The map $ \theta \colon ae \in Ae \mapsto ae\xi=a\xi \in F$
is onto because $\xi$ is a cyclic vector.
There are two choices for the $i$th restriction
mapping $\theta \restriction_{Ae_i}$ from $Ae_i$ to $F$.
Either it is one-to-one, or zero.  This is becuase the kernel
is an $A$-invariant subspace of $Ae_i$, and $Ae_i$ is a minimal
left ideal of $A$.  Let $e_{i_1}, \dots e_{i_m}$ be the ones
for which the restriction  map is one-to-one.
Then $\theta$ restricts to an isomorphism of $A$-modules
$Ae_{i_1} + \cdots + Ae_{i_m} \cong F$.

Since each $Ae_{i_j}$ is contained in the $\cal K$-module
${\cal K}e_{i_j}$, the direct sum
$Ae_{i_1} + \cdots + Ae_{i_m}$ is contained in the
direct sum of
$\cal K$-modules
${\cal K}e_{i_1} + \cdots + {\cal K}e_{i_m}$, and
we have proved that $F$ is contained in a $\cal K$-module.
$\Box$

\vskip\baselineskip

\noindent
{\bf Remark 1.2.  Modular Annihilator Algebras.} \
The subalgebra $A$ of $\cal K$ is a {\it
modular annihilator algebra}
as defined in [Barn, 1968], [Barn, 1969]
or [Palm, 1994], Chapter 8, Theorem 8.4.5.
We do not require
the ideal of finite ranks in $A$
(i.e. the socle) be dense in $A$'s topology.
But the assumption
of spectral invariance in $\cal K$ is strong enough to ensure that
$A$ be a modular annihilator algebra.
(See for example [Barn, 1968], Theorem 4.2.)

\vskip\baselineskip
\vskip\baselineskip
\section{References}
\smallskip
\smallskip
\footnotesize

\noindent[{\bf Barn, 1968}]
\, B. A. Barnes, 
{\it On the existence of minimal ideals in a Banach algebra},
Trans. Amer. Math. Soc. {\bf 133}
(1968), 511--517.

\noindent[{\bf Barn, 1969}]
\, B. A. Barnes, 
{\it Subalgebras of modular annihilator algebras},
Proc. Camb. Phil. Soc. {\bf 66}
(1969), 5--12.

\noindent[{\bf Palm, 1994}]
\, T. W. Palmer, 
{\it Banach Algebras and the General
Theory of $\star$-algebras,
Volume I: Algebras and Banach Algebras},
Encyclopedia Math. Appl. {\bf 49},
Cambridge Univ. Press, Cambridge, 1994.

\noindent[{\bf Sch, 1992}]
\, L. B. Schweitzer, 
{\it A short proof that $M_n(A)$ is local if $A$ is local and Fr\'echet}, 
Internat. J. Math.  {\bf 3} no. 4 (1992), 581--589.

\noindent[{\bf Sch, 2015}]
\, L. B. Schweitzer, 
{\it Extending representations of dense subalgebras of $C^\star$-algebras,
and spectral invariance}, 
arXiv:1502.02786v2.

\vskip\baselineskip
\noindent{Web Page: \url{http://www.svpal.org/~lsch/Math/indexMath.html}.}

\end{document}